\input amstex
\mag=\magstep1
\documentstyle{amsppt}
\nologo
\NoBlackBoxes
\NoRunningHeads
\TagsOnRight

\define\inbox#1{$\boxed{\text{#1}}$}
\def\fp{\flushpar}

\define\EE{\Bbb E}
\define\CC{\Bbb C}
\define\RR{\Bbb R}

\define\GL{\roman{GL}}
\define\SO{\roman{SO}}
\define\even{\roman{even}}
\define \dalpha{{\dot \alpha}}
\define \dbeta{{\dot \beta}}
\define \dgamma{{\dot \gamma}}

\define \Cliff{{\roman{Cliff}}}
\define \CLIFF{{\roman{CLIFF}}}
\define \CLIFFC{{\roman{CLIFF}^{\Bbb C}}}
\define \SPIN{{\roman{SPIN}}}
\define \END{{\roman{END}}}

\define \Spin{{\roman{Spin}}}

\define \SEE{{S\hookrightarrow \EE^n}}

\define \fg{\frak g}
\define \bee{{\bold e}}
\define \bEE{{\bold E}}
\define \ee{{\roman e}}

\define \tr{{\roman{tr}}}
\define \Ker{{\roman{Ker}}}
\define \s{{\sigma}}
\define \g{{\gamma}}
\define \ep{{\epsilon}}
\define \bb{{\overline b}}
\define\Not{\not \!\!}

\define\Bat{\not \!\!}

\font\twptbf=cmbx12

\def\fp{\flushpar}
\define\tp#1{\negthinspace\ ^t#1}

\define\lr#1{^{\sssize\left(#1\right)}}

\pagewidth{16.5 truecm}
\pageheight{23.0 truecm}
\vcorrection{-0.0cm}
\hcorrection{-0.5cm}

\document
\baselineskip=13.2pt 

\centerline{\twptbf
 Generalized Weierstrass Relation for a Submanifold $S^k$
 in $\EE^n$}

\centerline{\twptbf
 Coming from Submanifold Dirac Operator}
\vskip 0.5 cm

\vskip 0.5 cm
\centerline{Shigeki Matsutani}
\centerline{8-21-1 Higashi-Linkan Sagamihara 228-0811 JAPAN}
\centerline{e-mail: RXB01142\@nifty.com}
\vskip 1.0  cm

\centerline{\twptbf \S 1. Introduction}

Using the submanifold quantum mechanical scheme \cite{dC, JK},
the restricted Dirac operator in a $k$-spin submanifold immersed in
euclidean space $\EE^n$ ($0<k<n$) is defined \cite{BJ, Mat1-10}.
We call it submanifold Dirac operator.
Then it is shown that
the zero modes of the Dirac operator express the local
properties of the submanifold, such as the
Frenet-Serret and generalized Weierstrass relations.
In other words
this article gives a representation of a
further generalized Weierstrass relation for
the submanifold.

Before we start to explain the idea of the submanifold quantum
mechanics, we will recall three facts.

\roster

\item
Let us consider a left-differential ring $\bold P$ and
its element $Q$
 over a manifold $M$.
 In the book of Bj{\"o}rk \cite{Remark 1.2.16 in Bj},
it is stated that
 to treat its adjoint right operator  is difficult.
 Its essential is as follows.
Assume that $M$ is Riemannian.
For smooth functions $f_1$ and $f_2$ whose support is compact,
we consider the following integral as a bilinear form
of $f_1$ and $f_2$ formally
$$
\int_M dvol (f_1 Q f_2). \tag 1-1
$$
What is a natural adjoint of $Q$? One might regard
an action to $f_1$  obtained by
 partial integral as its adjoint.
However there exists an obstacle because the measure
depends upon the local coordinate.
Hence {\it concept of
 adjoint operator  is very subtle}.

\item
In a quantum mechanical problem, we sometimes
encounter the situation that for an eigen function
(and thus its zero mode) $\psi$ of a differential
operator $P$,
$$
P \psi = E \psi, \tag 1-2
$$
is a vector of a representation space of a group $G$.
In the case, if one finds a solution
of $P\psi =0$, he obtains
a representation of the group $G$.

Further suppose that
 $P$ is decomposed by $P = P_1 + P_2$.
Let us consider a kernel of $P_2$, $\Ker P_2$, in a certain
function
space. If an element $\psi_1 \in \Ker P_2$ satisfies
$$
	P_1 \psi_1 =E \psi_1, \tag 1-3
$$
we also obtain a representation of the group $G$.

\item In the quantum mechanical problem over a manifold $M$,
there are typical two pairings for the function space over $M$
in general, {\it i.e.},
1) global pairing $<,>$ induced from L$^2$-norm,
such as (1-1) and 2) point-wise pairing
$\cdot$ which is connected with the probability density.
(Of course, the point-wise pairing can be regarded as
$<\circ, \delta(p)\times > $ in terms of the Dirac distribution
$\delta()$ at $p$ in $M$.)

\endroster

Even though the concept of adjoint operator is subtle,
we could define a  natural adjoint operator if the measure
is fixed by some reason, {\it e.g.}, Haar measure.
We note that the ordinary Lebesgue measure in a euclidean
space $\EE^n$
is a typical Haar measure of the translation group.
The quantum mechanics in $\EE^n$ is based on the measure
and concept of the adjoint operator plays essential roles.
In such a case, by fixing a typical measure and
L$^2$-type paring $<,>: \Omega^* \times\Omega \to \CC$,
for an operator $Q\in \bold P$ whose domain is $\Omega$,
  we can define a right-adjoint operator, $Ad(Q)$,
with the domain $ \Omega^*$ by
$$
       <f,Q g> = < f Ad(Q), g> , \quad\text{for }
       (f,g)\in \Omega^* \times\Omega.
 \tag 1-4
$$
Assume that there is an isomorphism $\varphi$
between domains $\Omega$  and $  \Omega^*$ as a vector space,
$\varphi: \Omega \to \Omega^*$.
Then we can define the left-adjoint operator $Q^*$ as
$Q^*f := \varphi^{-1}(\varphi(f) Ad(Q))$.
Triplet
$(\Omega^*\times \Omega, <,>, \varphi)$ becomes
 a preHilbert space $\Cal H$ by letting
the inner product $(,)_\varphi:\Omega\times\Omega \to \CC$ be
$(f,g)_\varphi:=<\varphi(f),g>$.
(Then $(P^* f, g)_\varphi=(f,P g)_\varphi$.)

Suppose that $Q\in \bold P$ is {\it self-adjoint},
{\it i.e.}, the domains  of $Q^*$ and $Q$ coincide
and $Q^*=Q$ over there.
Then we have the following properties:

\roster

\item The kernel of $Q$, $\Ker Q$, is isomorphic to
the $\Ker (Ad(Q))$ i.e.,
$$
	(\Ker (Q))^* :=
	\varphi(\Ker (Q)) = \Ker( Ad(Q)).
        \tag 1-5
$$
\item $( (\Ker Q)^*\times \Ker Q, <,>, \varphi)$ becomes
 a preHilbert space.

\endroster

The projection $\pi$ from $\Omega^*\times\Omega$ to
$(\Ker Q)^* \times \Ker Q$ is
commutative with $\varphi$, i.e.
$$
	\varphi \pi|_{\Omega}=\pi|_{\Omega^*} \varphi , \quad
       (\varphi(\pi|_{\Omega} f) =
\pi|_{\Omega^*} \varphi(f) \equiv \varphi(f) Ad(\pi|_{\Omega})) .
\tag 1-6
$$
For $\pi$ satisfying (1-6), we will
say that {\it $\pi$ consists with the inner product}.
In fact (1-6) means that ${\pi_{\Omega}}^*={\pi_{\Omega}}$
due to the relation
${\pi|_{\Omega}}^*f \equiv \varphi^{-1}(\varphi(f)
Ad({\pi|_{\Omega}}))
={\pi|_{\Omega}} f$.

Next assume that $P_2$ is not self-adjoint in the preHilbert space
$\Cal H\equiv (\Omega^*\times \Omega, <,>, \varphi)$.
Even for the case, if $P_2$ is a certain operator,
 we could construct a transformation $\eta_{\roman{sa}}$
of the preHilbert space and its operators,
and find a preHilbert space $\Cal H'$ satisfying the following
conditions,
\roster

\item There exists an isomorphism $\eta_{\roman{sa}}:
\Omega^* \times \Omega \to \tilde \Omega^* \times\tilde  \Omega$
as a vector space.

\item By defining a pairing
$<\circ,\times>_{P_2}:=<\eta_{\roman{sa}}|_{\Omega^*}\circ,
\eta_{\roman{sa}}|_{\Omega}\times>$,
and $\tilde \varphi:=
\eta_{\roman{sa}}|_{\Omega^*} \varphi
 \eta_{\roman{sa}}^{-1}|_{\tilde\Omega}$,
$\Cal H'\equiv (\tilde \Omega^* \times\tilde  \Omega,
(,), \tilde\varphi)$.

\item An operator $P$ for $\Omega$ is transformed as
$\eta_{\roman{sa}}|_{\Omega^*} P
 \eta_{\roman{sa}}^{-1}|_{\tilde\Omega}$

\item $P_2$ itself is a self-adjoint in $\Cal H'$.
\endroster
We call $\eta_{\roman{sa}}$ {\it self-adjointization}:
$\eta_{\roman{sa}} : \Cal H \to \Cal H'$.
Of course, the self-adjointization is not a unitary operation
for $P_2$.

As mentioned above, $\Ker (P_2)\subset \tilde \Omega$ also
becomes a preHilbert space
denoted by $\Cal H'':=((\Ker(P_2))^* \times
\Ker (P_2), (,), \tilde \varphi)$.
Letting the projection of $\Cal H'\to \Cal H''$ be denoted by
$\pi_{P_2}$,  we have a sequence,
$$
	\Cal H {\overset \eta_{\roman{sa}} \to\longrightarrow}
 	\Cal H' {\overset \pi_{P_2} \to\longrightarrow}
	\Cal H''. \tag 1-7
$$
This sequence is a key of the submanifold quantum mechanics.
Instead of considering $P\psi =E\psi$ in $\Cal H$, we might
search a solution of $f$ in $\eta_{\roman{sa}}(P) \psi_1 =E\psi_1$
in $\Cal H''$.

Let us explain the idea of the submanifold quantum mechanics.
We note that for a smooth $k$-submanifold $S$ embedded
 in the $n$-euclidean space
$\EE^n$ ($0<k<n$), we can find a natural adjoint operator for a
differential operator defined over $S$ by
fixing the induced metric of $S$ from $\EE^n$,
even though $S$ is a curved space.
For a Schr{\"o}dinger equation in $\EE^n$ with
the L$^2$-type Hilbert space $\Cal H$,
$$
	-\Delta \psi = E\psi, \tag 1-8
$$
we should regard the Laplace operator
$\Delta$ as a Casimir operator for the
Lie group with respect to the translation.
By considering $\Delta$ over a tubular neighborhood of $S$,
$\Delta$ includes the normal differential operator
$\partial_\perp$.
We regard the normal differential
$\partial_\perp$ as above $P_2$.
As $\partial_\perp$ is not self-adjoint in general,
we step the above sequence.
In the self-adjointization $\eta_{\roman{sa}}$,  we obtain an
extra potential in the differential equation.
By considering kernel of purely normal component of the
differential $\partial_\perp$
and restricting the its definition region of
$\eta_{\roman{sa}}(\Delta)$
at $S$, we define a differential operator,
$$
	\Delta_{S \hookrightarrow \EE^n}
         :=\eta_{\roman{sa}}(\Delta)|_{\Ker\partial_\perp}|_S.
         \tag 1-9
$$
Then it turns out to be
$$
	\Delta_{S \hookrightarrow \EE^n}=
	\Delta_{S} + U(\kappa_i), \tag 1-10
$$
where $\Delta_{S}$ is the Beltrami-Laplace operator on $S$
which exhibits the intrinsic properties of $S$ and $U(\kappa_i)$
is an invariant functional of  principal
 curvature $\kappa$'s of $S$ in $\EE^n$.

Due to the self-adjointness of $\partial_\perp$ in $\Cal H''$,
we are allowed to consider the Hilbert space $\Cal H''$ for
$\Delta_{S \hookrightarrow \EE^n}$ naturally.
The point-wise product in $\Ker(\partial_\perp)|_S$
also has meaning;
the probability density is well-defined  there.
Thus we can consider the submanifold Schr{\"o}dinger equation,
$$
	-\Delta_{S \hookrightarrow \EE^n}\psi = E\psi,
$$
as a quantum mechanical problem and a representation of
translational group.

For the case  a smooth surface $S$ embedded in
$\EE^3$, by letting
$K$ and $H$ denote the Gauss and mean curvature,
we obtain
$$
	\Delta_{S \hookrightarrow \EE^3}
         =
	\Delta_{S} +  H^2-K. \tag 1-11
$$
Hence it is expected that
$ \Delta_{S \hookrightarrow \EE^n}$ and its zero mode
exhibit the extrinsic properties, {\it e.g.}, umbilical points,
 of the submanifold.

The submanifold quantum mechanics was opened by Jensen
and Koppe  about
thirty years ago and rediscovered by da Costa [JK, dC];
even though they did not mention essentials of the submanifold
quantum mechanics as described above, they obtained
$\Delta_{S \hookrightarrow \EE^n}$.
We should note that
{\it as the above operation is local, our consideration can be
extended to an immersed submanifold $S$ in $\EE^n$.}

\vskip 1.0 cm

As I have been considering the Dirac operator version of above
quantum system for this decade \cite{Mat1-10}, which is
our main subject in this paper.
In the investigation of the Dirac operator,
we should recall the  fact that solutions $\{\Psi\}$
of the Dirac equation,
$$
\Not {\bold D}_{\EE^n} \Psi=0,
$$
 locally represents the spin group.
By letting $\cdot$ denote point-wise pairing,
$(\varphi_{pt}(\{\Psi\}) \times\{\Psi\}, \cdot, \varphi_{pt})$
for a certain map $\varphi_{pt}$
becomes a preHilbert space and for an appropriate $\gamma$ matrix
and solution $\Psi$, $\varphi_{pt}(\Psi) \gamma \Psi$
exhibits a section of
the $\SO(n)$ principal bundle $\SO_{\EE^n}(T\EE^n)$ over $\EE^n$.

Thus we apply the submanifold quantum mechanical scheme to
the Dirac operator over a $k$-spin submanifold $S$
immersed in $\EE^n$.
Then we have a representation
 of $\SO_{\EE^n}(T\EE^n)|_S$  at $S$ immersed in $\EE^n$,
which is the generalized Weierstrass relation.
Our main theorem is Theorem 3.15.

The organization of this article is as follows.
Section 2 devotes the preliminary on
the geometrical setting \cite{E} and conventions of
the Clifford module and its related objects \cite{BGV, Tas}.
In \S 3, we give an construction algorithm of the submanifold Dirac
equation and investigate its properties.
Section 4 gives its example.

\vskip 0.5 cm
\centerline{\twptbf Acknowledgment}
\vskip 0.5 cm

I thank  the organizers of the conference,
especially Prof. M. A. Guest
and Prof. Y. Ohnita for giving me a chance to talk there.
It is acknowledged that Prof. M. A. Guest
gave me helpful suggestions to write this article.
I am grateful to Prof. B.~G.~Konopelchenko,
 Prof. U. Pinkall and
Prof. F.~Pedit for their interest in this work,
Prof. T. Kori and Dr. Y. Homma for inviting me
their seminar and crucial discussions.
I thank Prof. S. Saito, Prof. K. Tamano and H. Mitsuhashi
for critical discussions and encouragements.

\vskip 1.0  cm

\centerline{\twptbf \S 2. Preliminary}

As we consider a restriction of the differential operator
and the function space over a manifold, we should deal with
it in the framework of the sheaf theory and thus
 we employ the following notations.

\proclaim{2.1 Conventions}\cite{Mal}
{\rm
 For a fiber bundle $A$ over a differential manifold $M$ and an
open
set $U \subset M$, let $\Gamma(U,A)$ denote
 a set of smooth
sections of the fiber bundle $A$ over $U$.
Further for a point $p$ in $M$,
let $\Gamma(p,A)$ denote a stalk at $p$ of a set of smooth
sections of the fiber bundle $A$.
}
\endproclaim

Further we use Einstein convention and
let $\CC$ ($\RR$) denote the complex (real) field.
For brevity, we use the notations
$\partial_{u^\mu}:=\partial/\partial u^\mu$ for a certain parameter
$u^\mu$.  For a real number $x \in [n, n+1)$ and an integer $n$,
[x] denotes $n$. Let $\CC_M$ ($\RR_{M}$) denote
a complex (real) line bundle over a manifold $M$.

In order to define the submanifold Dirac operator,
let us consider a smooth spin $k$-submanifold $S$
immersed in $n$ euclidean space $\EE^n$.
As mentioned in \S 3, the Dirac operator can be constructed locally.
Thus in order to simplify the argument,
we assume that $S$ satisfies several properties as follows.

\proclaim{2.2 Assumptions/Notations on the
 submanifold $S$ and $T_S$}{\rm

\roster

\item
Let $S$ be diffeomorphic to $\RR^k$ as a chart
$(s^1,\cdots,s^k)\in \RR^k$.

\item
Let a tubular neighborhood of $S$ be denoted by
 $T_S$, $\pi_{T_S}:T_S \to S$.

\item Let $T_S$
 have a foliation structure by
$(dq^\dalpha=0)_{\dalpha=k+1,\cdots,n}$, $T_S \approx
\RR^{k}\times \RR^{n-k}$
$\ni(s^1,\cdots,s^k,q^{k+1},\cdots,q^n)\equiv(u^1,\cdots,u^n)$.
(Let beginning of Greek indices $\alpha$, $\beta$ run
from $1$ to $k$ and those with dot $\dalpha$, $\dbeta$
run from $k+1$ to $n$. Let the middle of Greek indices
$\mu$, $\nu$ run from $1$ to $n$.)

\item Let each leaf of the foliation
parameterized by $q=(q^{k+1}, \cdots q^n)$ be denoted by
$S_q$; $S_{q=0}\equiv S$.

\item
Let  $\fg_{T_S}$, $\fg_{S_q}$ and $\fg_{S}$ denote the
 Riemannian metric of $T_s$, $S_q$ and $S$
induced from that of
$\EE^n$ respectively.

\item Let
$g_{T_S \ i,j}:=\fg_{T_S}(\partial_{x^i},\partial_{x^j})$,
$g_{T_S \ \mu,\nu}:=\fg_{T_S}(\partial_{u^\mu},\partial_{u^\nu})$,
$g_{T_S}:=\det g_{T_S \ \mu,\nu}$
and so on for $T_S$, $S$ and $\EE^n$.

\item
Let $(dq^\dalpha)_{\dalpha=k+1,\cdots,n}$ be an orthonormal base
$\fg_{S_q}(\partial_{q^\dalpha},\partial_{q^\dbeta})
=\delta_{\dalpha,\dbeta}$
and satisfy \break
$\fg_{T_S}(\partial_{q^\dalpha},\partial_{s^\alpha})=0$
for $\alpha=1,\cdots,k$ and $\dalpha=k+1,\cdots,n$.
In other words, we have
$$
	\fg_{T_S}=\fg_{S_q}\oplus
     \delta_{\dalpha\dbeta} d q^\dalpha\otimes d q^\dbeta. \tag 2-1
$$

\item At a point $p$ in $S_q$, let an orthonormal frame
of the cotangent space $T^* S_q$ be denoted by
$d \xi:=( d \xi^\mu):=(d \zeta^\alpha, d q^\dalpha)$.
\endroster

We call this parameterization $q$ satisfying
(2)-(7){\it canonical parameterization}.

}
\endproclaim

\proclaim{\fp 2.3 Notation}
{\rm For a point $p$ of $S$, let
{\it Weingarten map} be denoted by
$-\gamma_\dbeta: T_p S \to T_p \EE^n$;
for bases $\bee_\alpha$ of $T S$ and
$\tilde\bee_\dbeta \in  T S^\perp$
($T_p\EE^n=T_p S \oplus T_p S^\perp$),
$$
     \gamma_\dbeta(\bee_\alpha) :=
\partial_\alpha \tilde\bee_\dbeta
    =
    \gamma^\dalpha_{\ \dbeta\alpha}
  \tilde\bee_{\dalpha}+\gamma^\beta_{\
    \dbeta\alpha} \bee_{\beta}
   .\tag 2-2
$$
}
\endproclaim

It is not trivial whether a submanifold in $\EE^n$
 has the canonical parameterization $q$
or not in general. However it is not so difficult to
prove
that a local chart of the submanifold
has canonical parameterization. In the proof,
2.2 (7) requires some  arguments on the Weingarten map but by
tuning
the frame $\tilde\bee_\dbeta$ in (2-2) using $\SO(n-k)$ action,
we can construct it due to the following Proposition \cite{Mat10},
which guarantees the existence of $S$ satisfying
the assumptions in  2.2.

\proclaim{\fp 2.4 Proposition}
{\rm For a  base $\bee_\alpha$ of $T S$,  there is
an orthonormal frame $\bee_\dalpha
= \delta_{\dalpha\dbeta} d q^\dbeta\in
T S^\perp$ satisfying
$$
    \partial_\alpha \bee_\dbeta =
    \gamma^\beta_{\ \dbeta\alpha}
  \bee_{\beta}.\tag 2-3
$$
}
\endproclaim
In terms of the properties, we have the moving
frame and the metric
as follows.

\proclaim{\fp 2.5 Lemma}
 For the moving frame
$\bee_\dalpha= \delta_{\dalpha\dbeta} d q^\dbeta$
 in  Proposition 2.4, the moving frame
 $\bEE_\mu=E_{\ \mu}^i \partial_i$,
$(E^i_{\ \mu}=\partial_\mu x^i)$
 in $S_q$ is expressed by
$$
    E^i_{\ \alpha} = e^i_{\ \alpha} +
            q^\dalpha \gamma^\beta_{\
            \dalpha\alpha} e^i_{\ \beta},
            \quad  E^i_{\ \dalpha} = e^i_{\ \dalpha}.
 \tag 2-4
$$
\endproclaim
\demo{Proof}
As a point $\bold x\equiv(x^i)$ in $S_q$
is expressed by
$
     \bold x= \bold y + \bee_\dalpha q^\dalpha
$
using $\bold y:=\pi_{T_S} \bold x$, we obtain them.
\qed \enddemo

\vskip 0.5 cm
\proclaim{\fp 2.6 Corollary}
 Let $g_{T_S}:=\det_{n\times n}(g_{T_S,\mu.\nu})$,
$g_{S_q}:=\det_{k\times k}(g_{S_q,\alpha.\beta})$ and
$\fg_{S}:=\fg_{S_q}|_{q=0}$.

\roster
\item  $g_{T_S}\equiv g_{S_q}$.

\item $\fg_{S_q}$ is expressed as
$
    \fg_{S_q}=\fg_{S}+
    \fg^{(1)}_{S} q^\dalpha+\fg^{(2)}_{S} (q^\dalpha)^2  ,
$
and locally
$$
g_{T {S_q}}(\partial_\alpha,\partial_\beta)
  = g_{S \alpha\beta}+
    [\gamma_{\ \dalpha\alpha}^\gamma g_{S\gamma\beta}+
   g_{S\alpha\gamma}\gamma_{\ \dalpha\beta}^\gamma]q^\dalpha
    +[\gamma_{\ \dalpha\alpha}^\delta g_{S\delta\gamma}
     \gamma_{\ \dbeta\beta}^\gamma]q^\dalpha q^\dbeta .
     \tag 2-5
$$

\item When we factorize $g_{S_q}$ as
 $g_{S_q}=g_{S}\cdot \rho_{S_q}$,
  the factor $\rho_{S_q}$ is
given by
$$ \split
     \rho_{S_q}=&1+2 \tr_{k\times k}
    (\gamma^{\alpha}_{\ \dalpha\beta})q^\dalpha
      +\left[ 2 \tr_{k\times k}(\gamma^{\alpha}_{\ \dalpha\beta})
          \tr_{k\times k}(\gamma^{\alpha}_{\ \dbeta\beta})-
          \tr_{k\times k}(\gamma^{\delta}_{\ \dalpha\beta}
         \gamma^{\alpha}_{\ \dbeta\delta})
         \right] q^\dalpha q^\dbeta\\ &+
         \Cal O(q^\dalpha q^\dbeta q^\dgamma) ,
\endsplit\tag 2-6
$$
where $\Cal O$ is Landau symbol.
\endroster
\endproclaim

Further we will recall properties of Clifford algebra
and its related quantities,
and show our conventions of them.

\proclaim{2.7 Clifford Algebra and Spinor Representations}{\rm
\cite{BGV, Tas}

\roster

\item
Let $\CLIFF(\RR^n)$ denote the Clifford algebra for
 the vector space $\RR^n$ and let \break
$\CLIFFC(\RR^n):=\CLIFF(\RR^n)\otimes \CC$.
Let $\CLIFFC^\even(\RR^n)$ denote the subaglebra of
$\CLIFFC(\RR^n)$
consisting of even degrees of generators in $\CLIFFC(\RR^n)$.
Further let $\SPIN(\RR^n)$ denote the spin group for $\RR^n$.

\item
For the exterior algebra
 $\wedge \RR^n= \oplus_{j=0}^n \wedge^j \RR^n$,
there is an isomorphism  as a vector space, called
symbol map $\CLIFF(\RR^n)\to \wedge \RR^n$.
Let  its inverse be denoted by $\gamma$,
$$
    \gamma:\wedge \RR^n \to \CLIFF(\RR^n),  \tag 2-7
$$
 which is called gamma-matrix.

\item For even $n$ case, let $\Cliff(\RR^n)$ denote
a $\CLIFFC(\RR^n)$-module
whose endomorphism is isomorphic to $\CLIFFC(\RR^n)$,
which is a $2^{[n/2]}$ dimensional $\CC$-vector space.
For odd $n$  case, let $\Cliff(\RR^n)$ denote
 a $\CLIFFC(\RR^n)$-module whose
endomorphism is homomorphic to $\CLIFFC(\RR^n)$ as a $2^{[n/2]}$
dimensional $\CC$-vector space representation and
elements are invariant for the action of
$\g(e_1)\g( e_2) \cdots \g(e_n)$.
Here $e_1, e_2, \cdots, e_n$ are orthonormal base of $\RR^n$.

\endroster
\endproclaim

In order to simplify the argument, we will fix the
expressions of the $\gamma$-matrices and so on as follows.

\proclaim{2.8 Conventions} \rm
\roster
\item
We recall the fact
$
\CLIFFC(\RR^{n+2}) \approx \END(\CC^2) \otimes \CLIFFC(\RR^n),
$
where $\END(\CC^2)$ is the endomorphism of $\CC^2$ and
$\END(\CC^2)$ can be generated by the Pauli matrices:
$$
\sigma_0:=\pmatrix 1 & 0 \\ 0 & 1 \endpmatrix, \quad
\sigma_1:=\pmatrix 0 & 1 \\ 1 & 0 \endpmatrix, \quad
\sigma_2:=\pmatrix 0 & -\sqrt{-1} \\ \sqrt{-1} & 0 \endpmatrix,
 \quad
\sigma_3
:=\pmatrix 1 & 0 \\ 0 & -1 \endpmatrix.\tag 2-8
$$

\item
For even $n$ case, noting
$\CLIFFC(\RR^{n}) \approx \END(\CC^{[n/2]})
$, we use the conventions,
$$
\split
\g(e_{1})&=\s_1\otimes\s_1\otimes\cdots\cdots\otimes\s_1
\otimes\s_1
\otimes\s_1\otimes\s_1, \\
\g(e_{2})&=\s_1\otimes\s_1\otimes\cdots\cdots\otimes\s_1
\otimes\s_1\otimes\s_1
\otimes\s_2,\\
\g(e_{3})&=\s_1\otimes\s_1\otimes\cdots\cdots\otimes\s_1
\otimes\s_1\otimes\s_1
\otimes\s_3,\\
\g(e_{4})&=\s_1\otimes\s_1\otimes\cdots\cdots\otimes\s_1
\otimes\s_1\otimes\s_2
\otimes\s_0\\
\g(e_{5})&=\s_1\otimes\s_1\otimes\cdots\cdots\otimes\s_1
\otimes\s_1\otimes\s_3
\otimes\s_0\\
 &\cdots \cdots \\
\g(e_{n-1})&= \s_1\otimes\s_3\otimes\cdots\cdots\otimes
\s_0 \otimes\s_0
\otimes\s_0 \otimes\s_0,\\
\g(e_{n})&=   \s_2\otimes\s_0\otimes\cdots\cdots\otimes
\s_0 \otimes\s_0
\otimes\s_0 \otimes\s_0.
            \endsplit\tag 2-9
$$

\item For odd $n$ case, noting
$
\CLIFFC(\RR^{n}) \approx \END(\CC^{(n-1)/2})\oplus
\END(\CC^{(n-1)/2})
$,
we use the conventions,
$$
\split
\g(e_{1})&=\s_1\otimes\s_1\otimes\cdots\cdots\otimes\s_1
\otimes\s_1\otimes\s_1
\otimes\frak e, \\
\g(e_{2})&=\s_1\otimes\s_1\otimes\cdots\cdots\otimes\s_1
\otimes\s_1\otimes\s_2
\otimes1 ,\\
\g(e_{3})&=\s_1\otimes\s_1\otimes\cdots\cdots\otimes\s_1
\otimes\s_1\otimes\s_3
\otimes1,\\
 &\cdots \cdots \\
\g(e_{n-2})&= \s_1\otimes\s_3\otimes\cdots\cdots\otimes
\s_0 \otimes\s_0
\otimes\s_0 \otimes1,\\
\g(e_{n-1})&=   \s_2\otimes\s_0\otimes\cdots\cdots\otimes
\s_0 \otimes\s_0
\otimes\s_0 \otimes1.\\
\g(e_{n})&=   \s_3\otimes\s_0\otimes\cdots\cdots\otimes
\s_0 \otimes\s_0
\otimes\s_0 \otimes1,
            \endsplit\tag 2-10
$$
where $\frak e$ is defined as the generator of
$\CLIFFC(\RR)=\RR[\frak e]/(1-\frak e^2)$.

\item By introducing
 $b_+:=\pmatrix 1 \\ 0\endpmatrix$
and $b_-:=\pmatrix 0 \\ 1\endpmatrix$,
$\Cliff(\RR^n)$ is spanned by the bases
$$
\Xi_{\epsilon}=b_{\ep_1}\otimes b_{\ep_2}\otimes\cdots\cdots
            \otimes b_{\ep_{[n/2]-1}}\otimes b_{\ep_{[n/2]}},
\tag 2-11
$$
where $\epsilon=(\ep_1,\ep_2,\cdots,\ep_{[n/2]})$ and
$\ep_a=\pm$ $(a=1,\cdots,[n/2])$.  Similarly, $\Cliff(\RR^n)^*$
is spanned by
$$
\overline\Xi_{\epsilon}=\bb_{\ep_1}\otimes \bb_{\ep_2}
\otimes\cdots\cdots
        \otimes \bb_{\ep_{[n/2]-1}}\otimes \bb_{\ep_{[n/2]}},
\tag 2-12
$$
where $\bb_+:=(1,0)$ and $\bb_-:=(0,1)$.
By appropriately numbering,
let $\epsilon^{[c]}$ ($c=1,\cdots,2^{[n/2]})$
denotes each $\epsilon$.
The isomorphic map
$\varphi:\Cliff(\RR^n)\to\Cliff(\RR^n)^*$ is given by
$$
\varphi(\sum_{c=1}^{2^{[n/2]}}a_c\Xi_{\epsilon^{[c]}})
    =\sum_{c=1}^{2^{[n/2]}}
\overline a_c\overline\Xi_{\epsilon^{[c]}},
      \tag 2-13
$$
for $a_c \in \CC$ and its complex conjugate
$\overline a_c$.

\item
Defining
 $b_1:=\dfrac{1}{\sqrt{2}}\pmatrix 1 \\ 1\endpmatrix$,
 $b_2:=\dfrac{1}{\sqrt{2}}\pmatrix 1 \\ \sqrt{-1}\endpmatrix$,
 $b_3:=\pmatrix 1 \\ 0\endpmatrix$,
$\bb_1:=\dfrac{1}{\sqrt{2}}(1,1)$, $\bb_2:=
\dfrac{1}{\sqrt{2}}(1,-\sqrt{-1})$
and $\bb_3:=(1,0)$,
$$
          \bb_a \sigma_b b_a=\delta_{a,b},
\quad(a,b=1,2,3), \ (\text{ not summed over }a). \tag 2-14
$$

\item
For $j=1,\cdots,[n/2]$, define
$$
     \Psi^{(2j-1)}:=b_1\otimes b_1\otimes \cdots \otimes
      b_1\otimes
      b_3\otimes b_1\otimes\cdots\otimes b_1,
$$
$$
     \Psi^{(2j)}:=b_1\otimes b_1\otimes \cdots \otimes b_1
        \otimes b_2
                  \otimes b_1\otimes\cdots\otimes b_1.
          \tag 2-15
$$
Further for odd $n$ case, we define,
$$
     \Psi^{(n)}:=b_3\otimes b_1\otimes \cdots \otimes b_1\otimes
b_1\otimes\cdots\otimes b_1,
            \tag 2-16
$$
and assume that $\frak e\Psi^{(i)}=\Psi^{(i)}$.
Further  define $\overline \Psi^{(k)}$ as in (2-13).
Then $\Psi^{(a)}$ is an element of $\Cliff(\RR^n)$
$(a=1,\cdots,n)$
satisfying
$$
\overline\Psi_{\{dx\}}^{(a)}\gamma(e^b)
 \Psi_{\{dx\}}^{(a)}=\delta^a_b,
\quad \text{ not summed over }a.\tag 2-17
$$

\endroster
\endproclaim

As we wish to consider a spin principal subbundle over $S$
induced from spin principal bundle over $\EE^n$,
we recall the facts:

\proclaim{\fp 2.9 Lemma }

\roster
\item
For $k < n$, $\CLIFFC(\RR^k)$ is a natural sub-vector
space of $\CLIFFC(\RR^n)$
by the inclusion of generators such that
$\CLIFFC^\even(\RR^k)$ is a subring of $\CLIFFC^\even(\RR^n)$.

\item
For $k < n$, $\SPIN(\RR^k)$ is a natural subgroup of
 $\SPIN(\RR^n)$.
\endroster
\endproclaim

\demo{Proof}
The statements are proved by considering five cases.
1)  $k=2l$, $n=2l+1$ case, 2) $k=2l$, $n=2l+2$ case,
3) $k=2l+1$, $n=2l+2$ cases, 4) $k=2l+1$, $n=2l+3$ case, and
5) otherwise. The fifth case can be proved by combinations
of the other cases. The first case is trivial.
The second case is a key because the third and forth cases
are similarly proved as the second case.
Therefore we concentrate our attention only on the
 $k=2l$, $n=2l+2$ case.
Recalling the facts in 2.8,
we have a natural inclusion as a set by the generators
for the bases $e_i$'s of $\RR^k$ and
$E_i$'s of $\RR^n$,
$$
        \tau_{k,n}: \g(e_i) \mapsto \g(E_i)
              := \sigma_1\otimes \g(e_i),
            \quad i=2,3, \cdots, k,  \tag 2-18
$$
and define $\tau_{k,n}( \g(e_i) \g(e_j) )
      := \tau_{k,n}( \g(e_i))\tau_{k,n}( \g(e_j) ) $,
 $\tau_{k,n}( \g(e_i) \g(e_j)\g(e_k) ) $ $
:= \tau_{k,n}( \g(e_i))$ $\tau_{k,n}( \g(e_j) )$
 $\tau_{k,n}(\g(e_k))$
and so on.
On the other hand, for $c^i\in \CLIFF(\RR^k)$, we have
a natural inclusion as an algebra,
$$
 \iota_{k,n}: \CLIFF(\RR^k)\hookrightarrow \CLIFF(\RR^n),
\quad   C_i = \sigma_0 \otimes c_i,\tag 2-19
$$
and then we have homomorphism,
 $\iota_{k,n}(c_i c_j)=\iota_{k,n}(c_i)\iota_{k,n}(c_j)$.
We note that for $1\le i,j\le k$,
$\g(E_i)\g( E_j)=\iota_{k,n}(\g(e_i) \g(e_j))$
and thus
in even subring $\CLIFF^{\even}(\RR^n)$,
image $\tau$ and $\iota$ agree. Thus (1) is proved.
 Accordingly
${\exp(\tau_{k,n}(\g(e_j)\g(e_i)))}$ can be regarded as an
elements of $\SPIN(\RR^n)$ and (2) is proved.
\qed \enddemo

\centerline{\twptbf \S 3. Construction of
Submanifold Dirac Operator}

An algorithm to construct submanifold Dirac operator
is the following six steps.

\subheading{ Step 1:
	Set the Dirac equation $\Not \bold D_{\EE^n}
        \Psi_{\EE^n} =0$
         in a euclidean space $\EE^n$
          and embed the $k$-smooth Spin submanifold $S$
          into $\EE^n$
           $(0<k<n)$}

Here we will give our notations in order to express the
Dirac equation $\Not \bold D_{\EE^n} \Psi_{\EE^n} =0$.

\proclaim{\fp 3.1 Clifford module etc.}\rm
\roster

\item Let $\CLIFFC_{\EE^n}(T^* \EE^n)$ denote the
Clifford bundle over $\EE^n$
which has the Clifford ring $\CLIFFC(\RR^n)$ structure
associated with
the cotangent bundle $T^*\EE^n$ and the set of
differential forms
$\Omega(\EE^n):=\sum_{a=1}^n \Omega^a(\EE^n)$.

\item Let $\Cliff_{\EE^n}(T^* \EE^n)$ denote the Clifford
module over $\EE^n$ associated with the
cotangent space $T^* \EE^n$, which is
 modeled by $\Cliff(\RR^n)$.

\item Let $\varphi_{pt}$ denote
the natural bijection between the spaces
as a $\CLIFFC_{\EE^n}(T^* \EE^n)$-module modeled by
(2-13), {\it i.e.},
$$
\varphi_{pt}:\Gamma(p, \Cliff_{\EE^n}(T^* \EE^n))\to
       \Gamma(p, {\Cliff_{\EE^n}(T^* \EE^n)}^*), \tag 3-1
$$
for a point $p$ in $\EE^n$.

\item
For an orthonormal frame
 $\bee^i\in T^*\EE^n$ at
$p\in \EE^n$,
let $\gamma_{\{\bee\}}$ denote
the $\gamma$-matrix as the map from $\Omega(\EE^n)$
to $\CLIFFC_{\EE^n}(T^* \EE^n)$ as (2-7),
$$
        \gamma_{\{\bee\}}: \bee^i \mapsto
        \gamma_{\{\bee\}}(\bee^i)
                 \in \CLIFFC_{\EE^n}(T^* \EE^n)|_p.\tag 3-2
$$
For later convenience, we also employ a notation for a one-from
$d u^\mu = E^\mu_i \bee^i\in \Gamma(p,\Omega^1(\EE^n))$,
$$
        \gamma_{\{\bee\}}( d u^\mu )
=E^\mu_i\gamma_{\{\bee\}} (\bee^i). \tag 3-3
$$
For simplicity, for a Cartesian coordinate system
$(x^i)=(x^1,\cdots x^n)$ in
$\EE^n$  we fix $\bee^i = d x^i$ and
in $T_S$, let $\bee^\mu = d \xi^\mu$ in 2.2 (8).

\item Let $\SPIN_{\EE^n}(T^* \EE^n)$ denote a spin principal
      bundle over $\EE^n$.
    Let a natural bundle map from $\SPIN_{\EE^n}(T^* \EE^n)$
      to $\SO(n)$-principal bundle $\SO_{\EE^n}(T^* \EE^n)$
       be denoted
      by  $\tau_{\EE^n}$; for $p \in\EE^n$,
      the orthonormal frame $\bee\in T^*_p\EE^n$ and
      $\ee^{\Omega}\in
      \SPIN_p(T^* \EE^n)$, the action of $\SO(n)$ is defined by
      $$
      \tau_{\EE^n}(\ee^{\Omega}) \ee^i
      :=\gamma_{\{\bee\}}^{-1}(
            \ee^{\Omega}\gamma_{\{\bee\}}(\ee^i)\ee^{-\Omega}).
             \tag 3-4
       $$

\item
The Dirac operator $\Bat {\bold D}_{\EE^n}$
 in the euclidean space $\EE^n$,
as an endomorphism between germs of
Clifford module $\Cliff_{\EE^n}(T^*\EE^n)$ over $\EE^n$,
$$
        \Bat {\bold D}_{\EE^n}:\Gamma(p,\Cliff_{\EE^n}(T^*\EE^n))
             \to \Gamma(p,\Cliff_{\EE^n}(T^*\EE^n)), \tag 3-5
$$
whose representation element is given by
$$
 \Not \bold D_{x,\{dx\}} = \gamma_{\{dx\}}( d x^i) \partial_{x^i} .
$$
\endroster\endproclaim

It is obvious that the following Proposition holds.

\proclaim{3.2 Proposition}
For a point $p\in \EE^n$ and a doublet
$(\overline \Psi_{\{dx\}},\Psi_{\{dx\}})$ in
$\Gamma(p, {\Cliff_{\EE^n}(T^* \EE^n)}^*)
        \times
\Gamma(p, \Cliff_{\EE^n}(T^* \EE^n))$,
the followings holds:
\roster

\item
$\delta_{ij}\overline \Psi_{\{dx\}}(x)
 \gamma_{\{dx\}}( d x^i)  \Psi_{\{dx\}}(x) d x^j
$ is an element of $\Gamma(p, \Omega^1(T^*\EE^n))$.

\item $\overline \Psi_{\{dx\}}(x)
  \Psi_{\{dx\}}(x)$
is an element of
 $\Gamma(p, \CC_{\EE^n})$.
\endroster
\endproclaim

Using the conventions (2-11),
the following proposition is easily obtained.

\proclaim{3.3 Proposition}
For a point $p$ in $\EE^n$,
$\{\Psi_{\{dx\}}^{[a]}
:=\Xi_{\epsilon^{[a]}}\}_{a=1,\cdots,2^{[n/2]}}$
satisfy the Dirac equation,
$$
\Not \bold D_{x,\{dx\}} \Psi_{\EE^n}(x) =0,
$$
and are bases of $\Gamma(p, \Cliff_{\EE^n}(T^* \EE^n))$ as
a $\Gamma(p,\CC_{\EE^n})$-vector space.
It means  that they are also the bases of the fiber
$\Cliff_{p}(T^* \EE^n)$ of $\Cliff_{\EE^n}(T^* \EE^n)$ at $p$.
\endproclaim

Further noting (2-17), we have the following important proposition.
\proclaim{\fp  3.4 Proposition}
There exist  $\Psi_{\{dx\}}^{(i)}$ $(i=1,\cdots,n)$ and
$\overline\Psi_{\{dx\}}^{(i)}:=\varphi_{pt}(\Psi_{\{dx\}}^{(i)})$
satisfying
$$
  \sum_{j,k=1}^n \overline\Psi_{\{dx\}}^{(i)}\gamma_{\{dx\}}(d x^j)
          \Psi_{\{dx\}}^{(i)}
           d x^k = d x^i.
$$
\endproclaim

We note that since
$\Cal W_p:=
       \Gamma(p, {\Cliff_{\EE^n}(T^* \EE^n)}^*)
        \times
\Gamma(p, \Cliff_{\EE^n}(T^* \EE^n))$
 has a point-wise-pairing $\cdot$
as mentioned in Proposition 3.2 (2), for a point $p\in \EE^n$,
$\Cal H_p:=(\Cal W_p,\cdot,\varphi_{pt})$
becomes the preHilbert space and due to Proposition 3.2 (1),
$\Cal W_p$ gives stalk $\Gamma(p,\SO_{\EE^n}(T^*\EE^n))$.

\subheading{ Step 2: By setting the canonical parameterization in
the tubular neighborhood $T_S$ of $S$, let the Dirac operator,
the Clifford module and so on be expressed by the parameterization}

Let
$\Gamma_{vc}(T_S,\Cliff_{\EE^n}(T^*\EE^n))$ denote a set of
sections of the
Clifford module $\Cliff_{\EE^n}(T^*\EE^n)$
whose support is in $T_S$. We go on to express its stalk at
$p\in T_S$ by $\Gamma_{vc}(p,\Cliff_{\EE^n}(T^*\EE^n))$.
Let $\Psi_{\{d \xi\}}$ denote a germ of
$\Gamma(p, \Cliff_{\EE^n}(T^*\EE^n))$ for a point $p \in T_S$.
Recalling (3-3), $dx^i=E^i_{\ \mu} du^\mu$ and decomposing
$E^{i}_{\ \mu} = G^{\ \nu}_\mu \Lambda^{i}_{\ \nu}$ as
$d x^i = \Lambda^i_{\ \mu} d \xi^\mu$,
we can find $\ee^\Omega\in \Gamma(p,\SPIN_{\EE^n}(T^*\EE^n))$
satisfying the following relations, called
{\it gauge transformation},
$$
        \Psi_{\{d \xi\}}(u)=\ee^{-\Omega} \Psi_{\{dx\}}(x), \quad
        \overline\Psi_{\{d \xi\}}(u)
        = \overline\Psi_{\{dx\}}(x)\ee^{\Omega},
$$
$$
\split
       \ee^{-\Omega} \gamma_{\{dx\}}( d x^i )
       \ee^{\Omega}=&
        \gamma_{\{d\xi\}}( {\Lambda}^{i}_{\ \mu}  d \xi^\mu )\\
      =& {\Lambda}^{i}_{\ \rho} G_\mu^{\ \rho} G_\nu^{\ \mu}
            \gamma_{\{d\xi\}}( d \xi^\nu )\\
      =&{E}_{\ \mu}^{ i} \gamma_{\{d\xi\}}( d u^\mu).
\endsplit \tag 3-6
$$
As $\overline\Psi_{\{d \xi\}}\equiv\varphi_{pt}(\Psi_{\{d \xi\}})$,
$\varphi_{pt}$ does not depend upon the orthonormal frame.
In terms of these expressions,
we have an assertion of Proposition of 3.2,
$$
\delta_{ij}\overline \Psi_{\{dx\}}(x)
 \gamma_{\{dx\}}( d x^i)  \Psi_{\{dx\}}(x) d x^j
=g_{T_S,\mu\nu}\overline \Psi_{\{d\xi\}}(u)
 \gamma_{\{d\xi\}}( d u^\mu)  \Psi_{\{d\xi\}}(u) d u^\nu,
$$
$$\overline \Psi_{\{dx\}}(x)
  \Psi_{\{dx\}}(x)=\overline \Psi_{\{d\xi\}}(u)
  \Psi_{\{d\xi\}}(u).\tag 3-7
$$
Further the representation of the Dirac operator
$\Not {\bold D}_{x,\{dx\}}$ is
transformed to $\Not \bold D_{u,\{d \xi\}}$ by means of
the gauge transformation,
$$
\split
	\Not \bold D_{u,\{d \xi\}}&= \ee^{-\Omega}
       \Not \bold D_{x,\{dx\}}
        \ee^{\Omega}
        =\gamma_{\{d \xi\}}( d u^\mu) \ee^{-\Omega}\partial_\mu
        \ee^{\Omega}\\
        &=\gamma_{\{d \xi\}}( d u^\mu)
           (\partial_\mu + \partial_\mu \Omega).
\endsplit \tag 3-8
$$
Then we have the following lemma:

\proclaim{3.5 Lemma}\it
We can regard $\Not \bold D_{u,\{d \xi\}}$ as a
representation of a map $\Bat \bold D_{T_S}$,
$$
        \Bat {\bold D}_{T_S}:\Gamma_{vc}(p,\Cliff_{\EE^n}(T^*\EE^n))
             \to \Gamma_{vc}(p,\Cliff_{\EE^n}(T^*\EE^n)),
$$
for a point $p$ in $T_S$

\endproclaim

Noting Lemma 3.5 and the fact that zero-section is in
$\Cliff_{\EE^n}(T^*\EE^n)$, it is not difficult to prove that
the solution space of
$\Not \bold D_{u,\{d \xi\}}\Psi=0$ in
$\Gamma(p,\CC^{2^{[n/2]}}_{\EE^n})$ belongs to
$\Gamma(p,\Cliff_{\EE^n}(T^*\EE^n))$.
Since $\Not \bold D_{u,\{d\xi\}}$ is a
$2^{[n/2]} \times 2^{[n/2]}$-matrix type first order differential
operator of rank $2^{[n/2]}$, the solution space
gives a $2^{[n/2]}$-dimensional
orthonormal frame in $\Gamma(p,\Cliff_{\EE^n}(T^*\EE^n))$
as a $\Gamma(p,\CC_{\EE^n})$-vector space.
In fact, due to Proposition 3.2, we can find one as
$\ee^{-\Omega} \Psi_{\{dx\}}^{[a]}$ by the gauge transformation.

Noting the Propositions 3.2-3.5 and (3-6)-(3-8),
we have a key proposition:

\proclaim{\fp 3.6 Proposition}

\roster

\item
For a point $p\in T_S$, there exist
$(\Psi_{\{d\xi\}}^{[a]})_{a=1,\cdots,2^{[n/2]}}$
in $ \Gamma(p,\CC^{2^{[n/2]}}_{\EE^n})$  such that
they are orthonormal frame  in $\Gamma(p,\Cliff_{\EE^n}(T^*\EE^n))$
and satisfy the Dirac equation
$$
\Not \bold D_{u,\{d \xi\}} \Psi_{\{d\xi\}}(u) =0, \tag 3-9
$$
i.e.,
$$
   \overline\Psi_{\{d\xi\}}^{[a]}
       \Psi_{\{d\xi\}}^{[b]} =\delta_{a,b},
\tag 3-10
$$
for
$\overline\Psi_{\{d\xi\}}^{[a]}:=
\varphi_{pt}(\Psi_{\{d\xi\}}^{[a]})$.

\item For
 $\Psi_{\{dx\}}^{[a]}\in \Gamma(p, \Cliff_{\EE^n}(T^*\EE^n))$
 in (1)
and  $\Psi_{\{dx\}}^{[a]}\in
 \Gamma(p, \Cliff_{\EE^n}(T^*\EE^n))$
in Proposition 3.3,
there exists a  germ
$\ee^{\Omega}\in$ $ \Gamma(p,\SPIN_{\EE^n}(T^* \EE^n))$
as a $(2^{[n/2]}\times2^{[n/2]})$-matrix satisfying
$$
\Psi_{\{d x\}}^{[a]}=\ee^{\Omega}\Psi_{\{d\xi\}}^{[a]},
\quad(a=1,\cdots,2^{[n/2]}).
$$

\item For $\ee^{\Omega}\in$ $ \Gamma(p,\SPIN_{\EE^n}(T^* \EE^n))$
in (2), by letting
$\Psi_{\{d\xi\}}^{(i)}:=\ee^{-\Omega}\Psi_{\{dx\}}^{(i)}$
for $(i=1,\cdots,n)$ in Proposition 3.4,
we have the relation,
$$
	g_{T_S,\mu\nu}\overline \Psi_{\{d\xi\}}^{(i)}(u)
 \gamma_{\{d\xi\}}( d u^\mu)  \Psi_{\{d\xi\}}^{(i)}(u) d u^\nu=dx^i,
\quad(i=1,\cdots,n). \tag 3-11
$$
\endroster
\endproclaim

\demo{Proof}
The proof of this proposition is done by the gauge transformation
except well-definedness of (3-11). We should check the dependence
of the orthonormal frame $\Psi_{\{d\xi\}}^{[a]}$.
Let us take another one $\Psi_{\{d\xi\}}^{\prime[a]}$.
By letting  $\ee^{\Omega'}$ defined as
$\Psi_{\{d x\}}^{[a]}=\ee^{\Omega'}\Psi_{\{d\xi\}}^{\prime[a]}$,
we have
$\Psi_{\{d \xi\}}^{[a]}=\ee^{-\Omega}\ee^{\Omega'}
\Psi_{\{d\xi\}}^{\prime[a]}$.
When one rewrites (3-11) in terms of
the frame $\Psi_{\{d\xi\}}^{\prime[a]}$,
there appears $\ee^{-\Omega'}\ee^{\Omega}
 \gamma_{\{d\xi\}}( d u^\mu)  \ee^{-\Omega}\ee^{\Omega'}$.
However as both are solutions of the Dirac equation (3-9),
the gauge transformation for $\ee^{-\Omega}\ee^{\Omega'}$
must leave the Dirac operator invariant.
Hence the $\partial_\mu$ component gives that
$\ee^{-\Omega'}\ee^{\Omega}
 \gamma_{\{d\xi\}}( d u^\mu)  \ee^{-\Omega}\ee^{\Omega'}$
$= \gamma_{\{d\xi\}}( d u^\mu)$.
 (3-11) does not
depend on the choice of the orthonormal frame. \qed
\enddemo

Here in order to investigate the domain of the new Dirac operator,
we wish to apply the facts in Lemma 2.9  to the fiber bundles.
Let
$\CLIFFC_{S}(T^*{S})$ (or $\Cliff_{S}(T^*{S})$)
be restricted its base space to the  submanifold $S$ by
$\CLIFFC_{\EE^n}(T^*\EE^n)|_{S}$
(or $\Cliff_{\EE^n}(T^*\EE^n)|_{S}$)
as a vector bundle over ${S}$. As the fiber of $\SPIN_p(T^*{S})$
is a subgroup of $\SPIN_p(T^*\EE^n)$, we also express
$\SPIN_{S}(T^*\EE^n) := \SPIN_{\EE^n}(T^*\EE^n)|_{S}$.
Here the reader should note the difference between
$\Gamma(p, \Cliff_{S}(T^*\EE^n))$ and
$\Gamma(p,\Cliff_{\EE^n}(T^*\EE^n))$; the former is
$\{\psi(s)\}$ and the later is $\{\psi(u)\}$.

\proclaim{\fp 3.7 Definition}
\roster

\item Let us denote the inclusion as a set due to the correspondence
of generators modeled by Lemma 2.9,
$$
  \tau_{{S},\EE^n}: \CLIFFC_{S}(T^*{S})  \to  \CLIFFC_{S}(T^*\EE^n),
\quad
(\tau_{{S},\EE^n}:\gamma_{{S},\{d \zeta\}}(d \zeta)
 \mapsto \gamma_{\{d \xi\}}(d \zeta)),
$$
so that this map induces a ring homomorphism from
$\CLIFFC_{S}^\even(T^*{S})  \to$\break
$  \CLIFFC^\even_{S}(T^*\EE^n)$.
Here $\gamma_{{S},\{d \zeta\}}$ denotes the $\gamma$-matrix
over $S$ associated with $T^* S$ and its orthonormal
 frame $\{d \zeta\}$.

\item Let us denote the  inclusion as a group modeled by Lemma 2.9,
$$
     \iota_{{S},\EE^n}: \SPIN_{S}(T^*{S})\to \SPIN_{S}(T^*\EE^n).
$$
\endroster\endproclaim

\subheading{ Step 3: In order to construct a
self-adjointization of the normal
differential operator, let us define a pairing in $T_S$}

As $q$ is a natural coordinate of $T_S$ and we are now
considering the affine geometry in category of differential
geometry, let us introduce the Haar measure with respect to
the local affine transformation along
the normal direction $\{q\}$, which is an invariant measure
for the rotation and translations of $q$, {\it i.e.},
$$
        \ee^{ b^\dalpha \partial_{q^\dalpha} }\
          g_{S}^{1/2}d^k s d^{n-k}q= g_{S}^{1/2}d^k s d^{n-k} q.
$$
Here $g_{S}^{1/2}$ is not a mistype of $g_{T_S}^{1/2}$.
We call this measure {\it normal affine invariance measure}.
Further we prepare a generic normal differential operator
$\partial_\perp:= b^\dalpha \partial_{q^\dalpha}$
for generic real constant number $b^\dalpha$'s.
As mentioned in Introduction, we apply the scheme in the
submanifold quantum mechanics to this operator.

\proclaim{3.8 Definition}\rm
For a point $p\in S$ and
$(\overline\Psi_{1, \{d\xi\}}, \Psi_{2, \{d\xi\}})
 \in \Gamma_{vc}(\pi_{T_S}^{-1}(p), \Cliff_{\EE^n}
 (T^*\EE^n)^*)$ $\times
\Gamma_{vc}(\pi_{T_S}^{-1}(p), \Cliff_{\EE^n}(T^*\EE^n))$,
we introduce L$^2$-type pairing $<|>$ in $T_S$ as a
fiber integral,
$$
<\Psi_{1, \{dx\}} |\Not \bold D_{u,\{d \xi\}}
\Psi_{2, \{d\xi\}}>
              = \int_{\pi_{T_S}^{-1}(p)} g^{1/2}_{T_S} d^{n-k} q
            \overline \Psi_{1, \{d\xi\}}(u)
 \Not \bold D_{u,\{d \xi\}}
\Psi_{2, \{d\xi\}}(u) \in \Gamma(p,\CC_S),
$$
and a transformation $\eta_{\roman{sa}}$ so that its
measure is the normal affine invariance measure,
$$
	\split
	<\Psi_{1, \{d\xi\}} |\Not \bold D_{u,\{d \xi\}}
          \Psi_{2, \{d\xi\}}>
              = &
	(\Phi_{1, \{d\xi\}} |\Not {\Bbb D}_{u,\{d \xi\}}
        \Phi_{2, \{dx\}})\\
              =& \int_{\pi_{T_S}^{-1}(p)} g^{1/2}_{S} d^{n-k} q
           \overline \Phi_{1, \{dx\}}(u)
              \Not {\Bbb D}_{u,\{d \xi\}}(u)
            \Phi_{2, \{dx\}}(u) \in \Gamma(p,\CC_S),
\endsplit
$$
where
$$
	\overline\Phi_{1, \{d\xi\}}
        :=\eta_{\roman{sa}}(\overline\Psi_{1, \{d\xi\}})
       \equiv
        {\rho_{S_q}}^{1/4}
           \overline\Psi_{1, \{d\xi\}},
\quad
	\Phi_{2, \{d\xi\}}:=\eta_{\roman{sa}}(\Psi_{2, \{d\xi\}})
       \equiv{\rho_{S_q}}^{1/4} \Psi_{2, \{d\xi\}},
$$ $$
           \Not \Bbb D_{u,\{d \xi\}}
:=\eta_{\roman{sa}}( \Not \bold D_{u,\{d \xi\}})
          \equiv {\rho_{S_q}}^{1/4}
              \Not \bold D_{u,\{d \xi\}}
             {\rho_{S_q}}^{-1/4}. \tag 3-12
$$
Further let $\tilde \varphi_{pt}$
denote $\eta_{\roman{sa}}\varphi_{pt}\eta_{\roman{sa}}^{-1}$.
\endproclaim

We should note that $\eta_{\roman{sa}}$ can be defined for
more general differential operators but for simplicity,
we only define it for the
Dirac operator. Further as the metric $\fg_{T_S}$ is not
singular, $\eta_{\roman{sa}}$ gives diffemorphism.
The following Lemma is naturally obtained.

\proclaim{3.9 Lemma } For points $p\in S_q$ and $p'\in S$,
the triplets $$
\Cal H_{p}':=
(\eta_{\roman{sa}}(\Gamma_{vc}(p, \Cliff_{\EE^n}(T^*\EE^n)^*)
\times
\Gamma_{vc}(p, \Cliff_{\EE^n}(T^*\EE^n)), \cdot,
\tilde \varphi_{pt}),
$$
and
$$
\Cal H':=(\eta_{\roman{sa}}
(\Gamma_{vc}(\pi_{T_S}^{-1}p', \Cliff_{\EE^n}(T^*\EE^n)^*)\times
\Gamma_{vc}(\pi_{T_S}^{-1}p', \Cliff_{\EE^n}(T^*\EE^n)), (|),
\tilde \varphi_{pt}),
$$
become preHilbert spaces.
In $\Cal H'$,
The operator $\partial_\perp$ is  anti-hermite, {\it i.e.},
$\partial_\perp^*=-\partial_\perp$ for each $b_\dalpha$.
\endproclaim

\subheading{ Step 4:
Decompose $\Not \Bbb D_{T_S}$ to normal part and tangential part}

Noting the orthonormal frame $d \xi$ in $T_S$ consisting of
$(d\zeta^1,\cdots,d\zeta^k,d q^{k+1},\cdots,d q^n)$,
let us decompose $\Not \Bbb D_{T_S}$ to
$$
       \Not \Bbb D_{u,\{d \xi\}}
      =\Not \Bbb D_{u,\{d \xi\}}^{\Vert}+
         \Not \Bbb D_{u,\{d \xi\}}^{\perp},
\tag 3-13
$$
where
$ \Not \Bbb D_{u,\{d \xi\}}^{\perp}
:= \gamma_{\{d\xi\}}( d q^\dalpha) \partial_{q^\dalpha}.
$
Then the following lemma is not difficult to proved.

\proclaim{3.10 Lemma}
\roster
\item $\Not \Bbb D_{u,\{d \xi\}}^{\Vert}$
 does not contain the vertical differential
operator  $\partial_{q^\dalpha}$.

\item ${\Not \Bbb D_{u,\{d \xi\}}^{\perp}}^*=-
 \Not \Bbb D_{u,\{d \xi\}}^{\perp}$ in $\Cal H'$.

\item For a point $p\in T_S$,
${\Not \Bbb D_{u,\{d \xi\}}^{\perp}}$
is a homomorphism of Clifford module,
$$
        {\Not \Bbb D_{u,\{d \xi\}}^{\perp}}:
        \eta_{\roman{sa}}(\Gamma_{vc}(p, \Cliff(T^* T_S))) \to
\eta_{\roman{sa}}(\Gamma_{vc}(p, \Cliff(T^* T_S))) .
$$

\endroster\endproclaim

For a point $p \in T_S$,
let us define a projection $\pi$,
$$
	\pi:\eta_{\roman{sa}}
          ((\Gamma_{vc}(p, \Cliff(T^* T_S)^*)\times
        \Gamma_{vc}(p, \Cliff(T^* T_S)))) \to
         ( \Ker_p(Ad({\partial_\perp}))\times
          \Ker_p(\partial_\perp)),
$$
where $\Ker_p(P)$ denotes the set of  germs of the kernel of $P$
at the point $p$ and $Ad(P)$ denotes the right-adjoint of $P$.
Here we note that $\Ker_p(\partial_\perp) (\subset$
$ \Gamma_{vc}(p, \Cliff(T^* T_S))$ is given as
the intersection  of $\Ker_p(\partial_\dalpha)$ for every
$\dalpha=k+1, \cdots, n$.

Then noting (2-6) and the fact that ${\rho_{S_q}}=1$ at $S$,
it is obvious that the following relations hold.
\proclaim{3.11 Lemma}

\roster
\item $\tilde \varphi_{pt}|_S =\varphi_{pt}|_S$
and $\eta_{\roman{sa}}(\Psi_{\{d\xi\}})|_S
\equiv \Psi_{\{d\xi\}}|_S $
for $p\in S$ and $\Psi_{\{d\xi\}}\in
\Gamma_{vc}(p,\Cliff_{\EE^n}(T^*\EE^n))$.

\item For a point $p\in T_S$,
$\tilde \varphi_{pt}|_{\Ker_p
(\partial_\perp)}$
is a bijection:
$$
   (\Ker_p(\partial_\perp))^*:=
   \tilde \varphi_{pt}\left(\Ker_p(\partial_\perp)
   \right)\approx \Ker_p(Ad(\partial_\perp)).
$$

\item $\pi$ consists with the inner product, $\pi^*=\pi$,
$$
	\pi \tilde\varphi_{pt}=\tilde\varphi_{pt} \pi,
         \quad \text{ and }\quad
	\pi \varphi_{pt}=\varphi_{pt} \pi \text{ at }S.
$$

\item For a point $p\in T_S$, the triplet
$
  \Cal H_p'':=((\Ker_p
      (\partial_\perp))^*\times
          \Ker_p(\partial_\perp), \cdot,
          \tilde\varphi_{pt}).
$
becomes a preHilbert space.
Let $\Cal H_p(\SEE):=\Cal H_p''$ for $p\in S$.

\item Any element in $\Ker_p(\partial_\perp)$
belongs to kernel of ${\not \Bbb D_{u,\{d \xi\}}^{\perp}}$,
{\it i.e.},
$$
\Ker_p(\partial_\perp)
\subset \Ker_p({\not \Bbb D_{u,\{d \xi\}}^{\perp}}).
$$
\endroster\endproclaim

\subheading{Step 5:
	Define submanifold Dirac operator ${\Not D}_{\SEE}$
at $\SEE$ by restricting its domain
$ \Ker_S(\partial_\perp)$}

Let us define
$\Not D_{\SEE}$ by restricting the its domain as
   ${\coprod_{p\in S}\Ker_p(\partial_\perp)}$,
{\it i.e.},
$$
	\Not D_{\SEE} := \Not \Bbb D_{T_S} |
   _{\coprod_{p\in S}\Ker_p(\partial_\perp)},
	 \tag 3-14
$$
with a local expression
$\Not D_{s,\{d \xi\}}:= \Not \Bbb D_{u,\{d \xi\}}
  |_{q=0, \partial_q=0}$
$\equiv\Not \Bbb D_{u,\{d \xi\}}^{\Vert}|_{q=0}$
whose explicit form is given by the following proposition.

\proclaim{\fp 3.12 Proposition }
For abbreviation of $\tau_{S,\EE^n}(\gamma_{S, d \zeta}(d
s^\alpha))$ by its image $\gamma_{\{d\xi\}}(d s^\alpha)$,
the explicit form of $\Not D_{s,\{d \xi\}}$ is given as
$$
\Not D_{s,\{d \xi\}}  =
          \tau_{S,\EE^n}(\Not {\bold D}_{S, s,\{d \zeta\}})
+         \gamma_{\{d \xi\}}( d q^\dalpha) \tr_{k\times k}
      (\gamma^\alpha_{\ \dalpha\beta}), \tag 3-15
$$
where $\Not {\bold D}_{S, s,\{d \zeta\}}$ is a proper
(or intrinsic)
Dirac operator of $S$,
$$
         \Not {\bold D}_{S, s,\{d \zeta\}}
        =
        \gamma_{\{d \xi\}}( d s^\alpha)(
        \partial_{s^\alpha} + \partial_{s^\alpha}\Omega|_{q=0}),
\tag 3-16
$$
by fixing the coordinate $s$ and the
orthonormal frame $\{d \zeta\}$.
\endproclaim

\demo{Proof}
\cite{BJ, Mat1-10, MT} From (2-6),
$
     {\rho_{S_q}}^{1/4}
        =1+\dfrac{1}{2} \tr_{k\times k}
    (\gamma^{\alpha}_{\ \dalpha\beta})q^\dalpha
          +   \Cal O(q^\dalpha q^\dbeta)
$.
Hence,
$$
{\rho_{S_q}}^{1/4}\partial_{\dalpha}{\rho_{S_q}}^{-1/4}
=\partial_{\dalpha}-
\frac{1}{2}\tr_{k\times k}    (\gamma^{\alpha}_{\ \dalpha\beta})
            + \Cal O(q^\dalpha),
$$
$$
{\rho_{S_q}}^{1/4}\partial_{\alpha}{\rho_{S_q}}^{-1/4}
=\partial_{\alpha}+\Cal O(q^\dalpha).
$$
The second term in (3-15) is obtained.
As $\tau_{S,\EE^n}$ induces the group inclusion $\iota_{S,\EE^n}$,
 we proved them. \qed \enddemo

Following proposition is important but is not difficult to be
proved.

\proclaim{3.13 Proposition}
For a point $S$, the stalks of
$\Gamma(p,\Cliff_{\EE^n}(T^*\EE^n))$
and $\Gamma_{vc}(p,\Cliff_{\EE^n}(T^*\EE^n))$
are bijective.
\endproclaim
As we are dealing with a germ at a point $p$  in $S$ hereafter,
let us neglect the difference between
$\Gamma(p,\Cliff_{\EE^n}(T^*\EE^n))$
and $\Gamma_{vc}(p,\Cliff_{\EE^n}(T^*\EE^n))$.

Let $p$ denote a point in $S$.
We are considering the kernel of
$\Not \Bbb D_{u,\{d \xi\}}^\perp$
at $p$ as a domain of the Dirac operator ${\Bat D}_{\SEE}$.
For given $\psi_0(s) \in \Gamma(p,\Cliff_S(T^*\EE^n))$,
we assume that $\Phi(u)\in \Gamma(p,\Cliff_{\EE^n}(T^*\EE^n))$
satisfies $\Not \Bbb D_{u,\{d \xi\}}^\perp\Phi=0$ with
a boundary condition $\Phi(s,0) =\psi_0(s)$ at $S$.
The existence of $\Phi$ is obvious because
we can find its solution as $\Phi(u)=\psi_0(s)$ due to
 $\partial_{q^\dalpha} \psi_0\equiv0$.
Therefore we regard that $\Gamma(p, \Cliff_S(T^*\EE^n))$
is a subset of $\Ker_p(\partial_\perp)$.
Further we can evaluate  $\Phi(u)\in
\Ker_p(\partial_\perp)$
around a point $p$ in $S$ as
$$
	\Phi(u)=\psi_0(s) + \sum \psi_\dalpha(s) q^\dalpha+
        \sum \psi_{\dalpha\dbeta}(s) q^\dalpha q^\dbeta + \cdots,
       \quad \psi_\dalpha(s)\equiv 0
         \tag 3-17
$$
where
$\psi_0(s) \in \Gamma(p,\Cliff_S(T^*\EE^n))$,
$\psi_\dbeta(s), \psi_{\dot \gamma \dot \delta}(s) \in
\Gamma(p,\CC_{\EE^n}^{2^{[n/2]}})$
and so on.
Accordingly we can identify
$\Ker_p(\partial_\perp)|_{q=0}$
with $\Gamma(p,\Cliff_S(T^*\EE^n))$ for a point $p\in S$.

Thus noting the fact
 that $\Not D_{s,\{d \xi\}}$ is expressed by coordinate-free
expressions and does not include $\partial_{q^\dalpha}$,
Definition 3.7, Lemma 2.9, 3.11 and Proposition 3.12
we have the following proposition:

\proclaim{3.14 Proposition}
The Dirac operator ${\Bat D}_{\SEE}$ is
 an endomorphism in germs of
Clifford module $\Cliff_{S}(T^*\EE^n)$ at a point $p$
in $S$,
$$
        {\Bat D}_{\SEE}:\Gamma(p,\Cliff_{S}(T^*\EE^n))
             \to \Gamma(p,\Cliff_{S}(T^*\EE^n)), \tag 3-18
$$
whose representation element is given by $\Not D_{s,\{d \xi\}} $.
\endproclaim

\subheading{Step 6:
	Consider solution of the
 submanifold Dirac equation $\Not D_{\SEE} \psi =0$
of $\SEE$}

Let us consider the solution of $\Not D_{\SEE} \psi =0$
at a point $p\in S$.  For non-vanishing $\psi \in
\Gamma(p,\CC_{\EE^n}^{2^{[n/2]}}))$ such that
$$
	\Not D_{s,\{d \xi\}} \psi(s) =0,
	 \tag 3-19
$$
$\psi(s)$ can be regarded as an element of
$\Gamma(S,\Cliff_{S}(T^*\EE^n))$.
By solving  a boundary problem for $\Phi(u) \in
\Gamma(S, \CC_{\EE^n}^{2^{[n/2]}})$,
$$
	\Not \Bbb D_{u,\{d \xi\}}^\perp\Phi(u) =0,\tag 3-20
$$
with boundary condition
$$
	\Phi|_S  = \psi, \text{ at } S, \tag 3-21
$$
 we have $\Phi(u)$ satisfying $\Not \Bbb D_{u,\{d\xi\}}\Phi(u) =0$
and $\Not \bold D_{u,\{d\xi\}}\eta_{\roman{sa}}^{-1}(\Phi(u)) =0$
at $p\in S$.

Noting that $\Not \Bbb D_{u,\{d\xi\}}|_{q\equiv0}
=\Not D_{s,\{d \xi\}}+\Not \Bbb D_{u,\{d\xi\}}^\perp$,
and $\Not D_{s,\{d \xi\}}$ does not include the parameter $q$,
we can apply the separation of variables to this system.
Thus it is expected that there exists an orthonormal frame
 $(\eta_{\roman{sa}}^{-1}
\Phi_{\{d\xi\}}^{[a]})_{a=1,\cdots,2^{[n/2]}}$ in
$\Gamma(p,\Cliff_{\EE^n}(T^*\EE^n))$ such that
each $\Phi_{\{d\xi\}}^{[a]}$ belongs to
$\Ker_p(\partial_\perp)$
and satisfies the Dirac equation (3-20);
each $(\eta_{\roman{sa}}^{-1}\Phi_{\{d\xi\}}^{[a]})$
is a solution of the Dirac equation (3-9).
Noting Lemma 3.11 (1),
we regard
$(\Phi_{\{d\xi\}}^{[a]}|_{q=0})_{a=1,\cdots,2^{[n/2]}}$
as the orthonormal frame of $\Gamma(p, \Cliff_{S}(T^*\EE^n))$
and solutions of the submanifold Dirac equation (3-19).

Inversely, as the Dirac operator $\Not D_{s,\{d \xi\}}$
is also a $2^{[n/2]}\times 2^{[n/2]}$-matrix type first order
differential operator, let us
assume that we find an orthonormal frame of
$\Gamma(p,\Cliff_S(\EE^n))$
 belonging to the solution space of the submanifold
Dirac equation (3-19).
 Let  it be denoted by $\psi_{\{d\xi\}}^{[a]}$
($a=1,\cdots,2^{[n/2]}$), {\it i.e.},
$\varphi_{pt}(\psi_{\{d\xi\}}^{[a]})
\psi_{\{d\xi\}}^{[b]}=\delta_{a,b}$.
For each $\psi_{\{d\xi\}}^{[a]}(s)$, we find an element
 $\tilde\Psi_{\{d\xi\}}^{[a]}(u)$ in the solution space of
 $\Not \Bbb D_{u,\{d\xi\}}\eta_{\roman{sa}}(\Psi(u)) =0$
or $\Not \bold D_{u,\{d\xi\}}\Psi(u) =0$ with
the boundary condition
$\tilde\Psi_{\{d\xi\}}^{[a]}(u)|_{q=0}=\psi_{\{d\xi\}}^{[a]}(s)$
at $p$.

Then due to Proposition 3.6 (2), we have an element of
$\ee^{\Omega} \in \Gamma(p, \SPIN_{\EE^n}(\EE^n))$ such that,
$$
	\Psi_{\{dx\}}^{[a]} =\ee^{\Omega}
           \tilde \Psi_{\{d\xi\}}^{[a]},
         \quad(a=1,\cdots,2^{[n/2]}).
$$
Using the element, we define $\Psi_{\{d\xi\}}^{(i)}
:=\ee^{-\Omega}\Psi_{\{d x\}}^{(i)}$,
$(i=1,\cdots,n)$. Then
Proposition 3.6 (3) gives the relation,
$$
 g_{T_S, \mu \nu} \overline\Psi_{\{d\xi\}}^{(i)}
        \gamma_{\{d \xi\}}(d u^\mu)
       \Psi_{\{d\xi\}}^{(i)} d u^\nu = d x^i ,
\quad(i=1,\cdots,n).\tag 3-22
$$
Therefore
by defining $\psi_{\{d\xi\}}^{(i)}(s)
 :=\Psi_{\{d\xi\}}^{(i)}(u)|_{q=0}$
and extracting the $d s^\alpha$ component,
we have
$$
	g_{S, \alpha\beta} \overline\psi_{\{d\xi\}}^{(i)}
        \gamma_{\{d \xi\}}(d{s^\alpha})
        \ee^{\Omega_i} \psi_{\{d\xi\}}^{(i)}
             = \frac{\partial x^i}{\partial s^\alpha},
\quad(i=1,\cdots,n,\alpha=1,\cdots,k).
         \tag 3-23
$$
This is a representation of the tangential space of $S$
or the generalized Weierstrass relation.
Here we note that $\psi_{\{d\xi\}}^{(i)}$ is also a solution of the
submanifold Dirac equation (3-19).
\vskip 1.0 cm

As mentioned above, our theory is locally constructed
and obtained objects do not depend on the local parameterization.
Accordingly we can easily extended it to that for a general
smooth spin $k$-submanifold in immersed $\EE^n$.
In other words, we regard the assumptions in 2.2 as
those for a local chart of the submanifold $S$.
Since $\EE^n$ is a spin manifold, we can define
 $\Cliff_{\EE^n}(T^*\EE^n)$,
$\SPIN_{\EE^n}(T^* \EE^n)$ and their global sections.
By restricting it on $S$, we have
$\Cliff_S(T^*\EE^n)$ and $\SPIN_S(T^*\EE^n)$ structure over $S$.
On the hand, $S$ has its own
$\Cliff_S(T^*S)$ and $\SPIN_S(T^*S)$ structure.
Due to Proposition 3.12,
$S$ has the submanifold Dirac operator $\Not D_{\SEE}$
with the local representation.

\proclaim{\fp 3.15 Theorem}
A smooth spin
$k$-submanifold $S$ immersed
in $\EE^n$ has a Clifford module
$\Cliff_S(T^*\EE^n)$ by locally constructing in terms of
 the solution space of the
submanifold Dirac equation $\Not D_{\SEE}\psi=0$ at each
chart and then the solutions give data of
tangential vector in $T\EE^n$ at each point as follows:

Let  $p\in S$ be  expressed by an affine coordinate $(x^i)$.
In a set $\{\psi\}$ of non-vanishing
germs of $\Gamma(p, \CC_S^{2^{[n/2]}})$
satisfying $\Not  D_{s,\{d\xi\}}\psi=0$,
there exists an orthonormal frame
$\{\psi_{\{d\xi\}}^{[a]}\}_{a=1,\cdots,2^{[n/2]}}$
$(\subset \{\psi\})$ of $\Gamma(p, \Cliff_S(T^*\EE^n)$
as a $\Gamma(p,\CC_{\EE^n})$-vector space.
For these bases, there exists a germ
$\ee^{\Omega} \in \Gamma(p, \Spin_{\EE^n}(T^* \EE^n))$
such that they are transformed to
$\{\Psi_{\{d x\}}^{[a]}|_S\}_{a=1,\cdots,2^{[n/2]}}$
by $\psi_{\{d\xi\}}^{[a]}
=\ee^{-\Omega}\Psi_{\{d x\}}^{[a]}|_S$,  $(a=1,\cdots,2^{[n/2]})$.
We define
$\psi^{(i)}:=\ee^{- \Omega} \Psi_{\{dx\}}^{(i)}|_S$ and
$\overline  \psi^{(i)}
:=\overline \Psi_{\{dx\}}^{(i)}\ee^{-\Omega} |_S$
$(i=1,\cdots,n)$.
Then the following relation holds:
$$
 g_{S, \alpha,\beta}\overline\psi^{(i)} [\tau_{S,\EE^n}
(\gamma_{S, \{d \xi\}}(d {s^\beta}))] \psi^{(i)}
           =\partial_{s^\alpha} x^i,
   \quad(i=1,\cdots,n, s=1,\cdots,k). \tag 3-24
$$

\endproclaim

This is the generalized
Weierstrass relation for a spin submanifold
immersed in $\EE^n$. Since the equation is given by the local
argument, the proof of the theorem  does not require a
global argument at all.

We note that even though we fix the coordinate system,
for a point $p$ in $S$,\break
$\left(\dfrac{\partial x^i}{\partial s^\alpha}\right)
_{\alpha=1,\cdots,k, i=1,\cdots, n}$ gives a data of
embedding $T^*_pS\approx \RR^k$ in $T^*_p\EE^n \approx \RR^n$
up to action of $\GL(\RR^k)$ and  $\GL(\RR^{n-k})$.
In other words, the submanifold Dirac operator
$\Not D_{\SEE}$ exhibits
data of Grassmannian bundle over $S$.

As we have assumed the canonical parameterization on
$T_S$ in the definition of the submanifold Dirac operator,
it should be extended to one without the assumption.
Then we can argue the global behavior of the submanifold Dirac
operator $\Not D_{\SEE}$ and its zero modes  over $S$.
It is expected that the zero modes bring us a
global data of the orthonormal frame of
$\Cliff_S(T^*\EE^n)$, {\it e.g.}, topological defect and so on.
Of course, for the low codimension cases, we can find the global
form of $\Not D_{\SEE}$ easily and then we can also argue it.
Inversely, if a $k$-spin manifold $M$ and a differential operator
over $M$ satisfy appropriate conditions, we also obtain the data
of immersion of $M$ in $\EE^n$ following the philosophy of
the Weierstrass relation.

\proclaim{\fp 3.16 Remark}\rm
The generalized Weierstrass relation for a conformal surface
 in $\EE^3$
was discovered by Kenmotsu [Ke] as a generalization of
 the Weierstrass relation for a minimal surface in  $\EE^3$.
The Dirac type relation was found by Konopelchenko
in 1995 [Ko1] and I showed that the submanifold Dirac equation
is identified with the relation \cite{Mat8}.
I also computed the submanifold Dirac operators for
the generalized Weierstrass relation for a conformal surface in
$\EE^4$,
which was also discovered by Konopelchenko [KL, Ko1-2]
and Pinkall
 and Pedit [PP].
Recently there are so many studies on the relations
between submanifolds and Clifford bundles
[Bo, Ko1-2, KL, KT, Fr, Tai1-2, Tr].
However no article mentioned the transformation $\eta_{\roman{sa}}$
as long as I know.
For this decade,
I have applied the transformation to the relation
between the Dirac operators and submanifolds, especially curves
immersed in
$\EE^n$. For the curve case, the submanifold Dirac equation is
identified with
Frenet-Serret relation \cite{Mat1-7}.
However  I could not reach the direct
connection between a more general
submanifold and zero modes of the Dirac operator
because the net meaning of the transformation
remained as a question in my mind. After the conference,
I noticed the essentials of the submanifold quantum mechanics
as mentioned in the Introduction,
which naturally leads me to the further generalized Weierstrass
relation or Theorem 3.15.

Though we did not mention in this article,
it is known, from the physical point of view,
that the obtained Dirac operator has the following properties;
\roster

\item The index of the Dirac operator
is related to the topological index for the case of curves
\cite{Mat2, 5}.

\item
The operator determinants are associated with the energy
functionals,
such as Euler-Bernoulli
functional and Willmore functional for the cases of
 space curves and of immersed
conformal surfaces respectively \cite{Mat6, 10}.

\item
The deformations preserving all eigenvalues of the submanifold
Dirac operators become the soliton equations,
such as MKdV equation \cite{Mat2, 6, 7, MT}, complex MKdV equation,
Nonlinear Schr{\"o}dinger equation \cite{Mat1, 3}, modified
Novikov-Veselov equations \cite{Ko1, Ko2, KL, KT, Tai1}
depending on submanifolds.

\endroster

As we showed in this article,
we modeled the theory of Thom class \cite{BT, BGV}
on construction of
this theory of the submanifold
Dirac operator. Both stand upon
vary similar geometrical situation.
As the Dirac operators are generally related to
some characteristic class of fiber bundle,
it is expected that the submanifold
Dirac operator might be related to
the Thom class and/or generalization of
Riemann-Roch [PP].

Further our scheme does not need global properties
and thus might be extended to a  subgroup
manifold immersed in more general group manifold
with a Casimir operator and a Haar measure
as  a continuous
group version of induced representation for finite group.

\endproclaim

\vskip 0.5 cm
\centerline{\bf \S 4. Dirac Operator on a conformal
 surface in $\EE^4$}
\vskip 0.5 cm

As an example, we will consider the case of a conformal surface
immersed in $\EE^4$.
In \cite{Mat10}, we gave an explicit local form (4-1) of
the submanifold Dirac
operator in this case and a conjecture
that the submanifold Dirac operator represents the surface.
 As the
conjecture was actually proved by Konopelchenko [KL, Ko2]
and Pedit and
Pinkall [PP], we will give another proof of the conjecture
 by means of
the submanifold Dirac system method. This means that my
conjecture in \cite{Mat10} was based upon physical assurance.

First we will give the properties of the Dirac operators
in a conformal
surface \cite{P}. For the case of a conformal surface $S$
in $\EE^n$, we can
set the metric given by
$$
        g_{S \alpha\beta}= \rho \delta_{\alpha \beta},
$$
and the orthonormal frame $\{d \zeta\}$ given by
$ d \zeta^\alpha:=
\rho^{-1/2} d s^\alpha$. Let complex parameterization of
$S$, $d z :=
d s^1 \pm\sqrt{-1} d s^2$.
We introduce another transformation
$\eta_{\roman{conf}}$ as follows. For a point $p\in S$ and
$\psi_{\{d \xi\}}\in \Gamma(p, \Cliff_S(T^*\EE^n))$,
let $\varphi_{\{d
\xi\}} := \eta_{\roman{conf}}({\psi}_{\{d \xi\}})
\equiv\rho^{1/2} \psi_{\{d \xi\}}$ and
 $\overline \varphi_{\{d\xi\}}
:= \eta_{\roman{conf}}(\overline{\psi}_{\{d \xi\}})
\equiv\varphi_{pt}( \psi_{\{d \xi\}})$.
 Then we have the
following properties.

\proclaim{\fp 4.1 Lemma}

\roster
\item The proper Dirac operator of $S$ is given by
$$
        \Not \bold D_{s,\{d\zeta\}}= \rho^{-1}
              \sigma^\alpha \partial_\alpha \rho^{1/2},
$$
by letting $\gamma_{S, \{d\zeta\}}(d\zeta^a)=\sigma^a$.

\item For a doublet
$(\overline \varphi_{\{d\xi\}},\varphi_{\{d\xi\}})$
in $\eta_{\roman{conf}}\Cal H_p(\SEE)$, $\varphi_{\{d\xi\}}
\tau_{S,\EE^n}[\gamma_{\{d \xi\}}
 ( d \xi^\alpha)] \varphi_{\{d\xi\}}
d s^\alpha$ is invariant form for choice of
 local parameterization of $S$.
\endroster\endproclaim

\demo{Proof}
(1) is obvious [P, Mat8-9]. By setting $\gamma_{\{d
\xi\}}(d s^\alpha)= \rho^{-1/2}\gamma_{\{d \xi\}}
(d \xi^\alpha)$, the
relation in Proposition 3.2 (1) becomes
$$
        g_{S \alpha \beta} \overline \psi_{\{d\xi\}}
           \tau_{S,\EE^n}[\gamma_{\{d \xi\}}
            ( d s^\alpha)] \psi_{\{d\xi\}}d s^\alpha
            =\overline \varphi_{\{d\xi\}}
        \tau_{S,\EE^n}[
         \gamma_{\{d \xi\}}(d \xi^\alpha) ]
        \varphi_{\{d\xi\}}d s^\alpha,
$$
and thus (2) is proved\qed \enddemo

\proclaim{\fp 4.2 Proposition}
\roster

\item For the case of a conformal surface in $\EE^4$,
 the submanifold
Dirac operator is given by
$$
\Not D_{s,\{d \xi\}} = 2\pmatrix \ & \ & ]
            \overline p_c & \partial  \\
                \ & \ & \overline\partial & - p_c\\
                       p_c & \partial & \ & \ \\
        \overline  \partial & -\overline p_c & \ & \
        \endpmatrix, \tag 4-1
$$
where $\partial:=(\partial_{s^1}-\sqrt{-1}\partial_{s^2})/2
$, $\overline\partial:=(\partial_{s^1}+\sqrt{-1}\partial_{s^2})/2$
and $p_c$ is
$$
        p_c:= -\frac{1}{2}\rho^{1/2}
             \tr_{2\times 2}(\gamma_{3\beta}^{\alpha}
                      +\sqrt{-1}\gamma_{4\beta}^{\alpha}).
$$

\item Let $\EE^4\approx \CC\times \CC\in (Z^1, Z^2)
\equiv (x^1 +\sqrt{-1} x^2, x^3 +\sqrt{-1} x^4)$. For the affine
coordinate $(Z^1,Z^2)$ of the surface, the relations,
$$
        d Z^1 = f m d z - g n d\overline   z,
        \quad
        d Z^2 =  f\overline  ndz + g\overline m d \overline z,
    \quad
        d\overline  Z^1 = \overline{d Z^1},\quad
        d\overline  Z^2 = \overline{d Z^2},
$$
hold if
$\varphi_1:=\pmatrix f\\ g\\ 0 \\ 0 \endpmatrix$
$\varphi_2:=\pmatrix  0 \\ 0 \\ m \\ n \endpmatrix$
are solutions of
$$
        \Not D_{s,\{d \xi\}} \varphi_a =0,
$$
and $(|f|^2+|g|^2)(|m|^2+|n|^2)=\rho^{1/2}$.
\endroster\endproclaim

Direct computation leads the next lemma.
\proclaim{\fp 4.3 Lemma }
For the four dimensional euclidean space $\EE^4$,
we can fix the representation of the gamma matrices
of $\{d x^i\}$ system as
$\gamma_{\{d x\}}(d x^i)= \sigma^1\otimes\sigma^i$
for $i=1,2,3$ and
$\gamma_{\{d x\}}(d x^4)= \sigma^2\otimes 1$.
By letting
$$
\Psi_1:=\pmatrix 1 \\0\\ 1\\ 0
\endpmatrix, \quad
\Psi_2:=\pmatrix 0 \\1\\ 0\\1 \endpmatrix, \quad
\Psi_3:=\pmatrix 1 \\0\\ 0\\1 \endpmatrix, \quad
\Psi_4:=\pmatrix 0 \\1\\ 1\\0\endpmatrix, \quad
$$
$\overline\Psi_1 = (0, 1, 0, 1)$, $
        \overline\Psi_2 = (1, 0, 1, 0)$,
        $\overline\Psi_3 = (0, -1, 1, 0 )$,
and      $ \overline\Psi_4 = (1, 0, 0, -1)$,
we have the following relations:
$$
 \sum_i \overline \Psi_1 \gamma_{\{d x\}}
      (d x^i) \Psi_1 d x^i = 2d Z_1, \quad
\sum_i
       \overline \Psi_2 \gamma_{\{d x\}}(d x^i) \Psi_2 d x^i
         = 2d \overline Z_1,
$$
$$\sum_i
       \overline \Psi_3 \gamma_{\{d x\}}(d x^i)
       \Psi_3 d x^i = 2d  Z_2,
\quad
\sum_i
       \overline \Psi_4 \gamma_{\{d x\}}(d x^i) \Psi_4 d x^i
        = 2d \overline Z_2.
$$

\endproclaim

\demo{Proof of Proposition 4.2}
(1) is proved by Theorem 3.15. We consider (2).
For $\varphi_a$'s, their independent partner solutions
are given by
$\varphi_3:=\pmatrix-\overline g\\
\overline f\\ 0 \\ 0 \endpmatrix$
$\varphi_4:=\pmatrix 0 \\ 0 \\
-\overline n \\ \overline m \endpmatrix$.
Since we have fixed
$\gamma_{S, \{d \xi\}}(d {\zeta^\beta})=\sigma^\beta$,
the convention in 2.8 gives
$$
\tau_{S,\EE^2}
(\gamma_{S, \{d \xi\}}(d {\zeta^\beta}))=
 \sigma^1 \otimes \sigma^\beta .
$$
We  define $\tilde \varphi_1:=\varphi_1+\varphi_2$,
$\tilde \varphi_2:=\varphi_3+\varphi_4$, $\tilde
\varphi_3:=\varphi_1+\varphi_4$,
and $\tilde \varphi_4:=\varphi_3+\varphi_2$.
Let us assume that $\tilde\varphi_a =
\rho^{1/2}\ee^\Omega \Psi_a|_{q=0}$
$(a=1,2,3,4)$ of Lemma 4.3. Then we find the spin matrix as
$$
  \rho^{1/2}\ee^\Omega=\pmatrix f & -\overline g& 0 & 0 \\
            g & \overline f& 0 & 0 \\
            0 & 0 & m & -\overline n \\
            0 & 0 & n & \overline m
        \endpmatrix.
$$
We have these dual bases,
$
        \overline {\tilde \varphi}_a =
        \overline \Psi_a \ee^{-\Omega}|_{q^\dbeta=0},
$ and obtain the relation,
$$
        2d Z_1 = \overline{\tilde \varphi}_1
         \sigma_1\otimes \sigma^\alpha
        {\tilde \varphi}_1 d s^\alpha,\quad
        2d \overline Z_1 =
         \overline{\tilde \varphi}_2\sigma_1\otimes
         \sigma^\alpha{\tilde \varphi}_2 d s^\alpha,
$$
$$
        2d Z_2 = \overline{\tilde \varphi}_3\sigma_1
            \otimes \sigma^\alpha
        {\tilde \varphi}_3 d s^\alpha,\quad
        2d \overline Z_2 = \overline{\tilde \varphi}_4
        \sigma_1\otimes
        \sigma^\alpha{\tilde \varphi}_4 d s^\alpha.
$$
Explicit representation of them proves (2).\qed
\enddemo

\enddocument

\bye


\Refs
\widestnumber\key{BBEIM}



\ref \key BGV \by N.~Berline, E.~Getzler and M.~Vergne
 \book Heat Kernels and Dirac Operators
\publ Springer \yr 1996 \publaddr Berlin \endref

\ref
   \key  Bj
   \by Bj\"ork, J-E
    \book Analytic $\Cal D$-Modules and Applications
    \publ Kluwer \publaddr Dordrecht, \yr 1992
\endref

\ref \key Bo \by A.~I.~Bobenko \paper
Surfaces in terms of 2 by 2 matrices: Old and new
 integrable cases
\inbook Harmonic Maps and Integrable Systems \eds
A.~P.~Fordy and J.~C.~Wood
\publ Vieweg \publaddr Wolfgang Nieger \yr 1994 \endref

\ref \key BJ \by M.~Burgess and B.~Jensen
  \jour Phys. Rev. A
\paper Fermions near two-dimensional sufraces
\vol48  \yr1993\page1861-1866 \endref


\ref \key  dC \by R. C. T.  da Costa \yr 1981
\jour Phys. Rev A \vol 23\pages 1982-7
\paper Quantum  mechanics of a constrained particle\endref

\ref \key E \by L.~P.~Eisenhart \book
A Treatise on the Differential Geometry
\publ Ellis Horwood \publaddr New York \yr 1909 \endref

\ref \key Fr \by T.~Friedrich \paper
On the Spinor Representation of Surfaces
in Euclidean 3-Space \jour J. Geom. Phys.
\vol 28 \pages 143-157 \yr 1997 \endref

\ref \key JK \by H.~Jensen and H.~Koppe
\paper Quantum Mechanics with Constraints
\jour Ann. Phys. \vol 63 \yr 1971 \pages 586-591 \endref

 \ref \key Ke \by K.~Kenmotsu \paper
Weierstrass formula for surfaces of prescrived mean curvature
\jour Math. Ann. \yr 1979 \pages89-99 \vol 245 \endref

\ref \key Ko1 \by B.~G.~Konopelchenko
 \jour Studies in Appl.~Math.
\paper Induced  surafces and their integrable dynamics
\vol 96  \yr1996 \pages 9-51 \endref

\ref \key Ko2 \bysame
\jour Ann. Global Analysis and Geom. \paper
Weierstrass representations for surfaces in 4D spaces
and their integrable deformations via DS hierarchy
\vol 18 \pages 61-74 \yr 2000 \endref

\ref \key KL \by B.~G.~Konopelchenko and G.~Landolfi
\jour Studies in Appl.~Math. \paper
Induced Surfaces and Their Integrable
Dynamics II.
Generalized Weierstrass Representations in 4D spaces
and Deformations via DS Hierarchy
\yr 2000 \vol 104 \pages 129-169 \endref

\ref \key KT \by B.~G.~Konopelchenko and I.~A.~Taimanov
 \jour J.~Phys.~A: Math.~\& Gen.
 \paper Constant mean curvature surfaces via
 an integrable dynamical system
\vol 29  \yr1996 \page1261-65 \endref

\ref \key Mal\by A.~Mallios \book Geometry of
Vector Sheaves, An Axiomatic
Approach to Differential Geometry I:
 Vector Sheaves. General Theory \publ
Kluwer   \yr 1998 \publaddr Netherlands \endref

\ref \key Mat1 \by S.~Matsutani
      \paper Berry phase of Dirac particle in thin rod
      \jour  J.  Phys.  Soc.  Jpn.  \vol 61 \yr 1992 \pages
        3825-3826 \endref



\ref \key Mat2 \bysame \paper The Relation
 between the Modified Korteweg-de
Vries Equation  and Anomaly of Dirac Field on a
Thin Elastic Rod
      \jour  Prog.Theor.  Phys.
       \vol 5 \yr 1994 \pages  1005-1037 \endref

\ref \key Mat3 \bysame
      \paper On the physical relation between the
Dirac equation and the
      generalized mKdV equation on
       a thin elastic rod
      \jour  Phys.  Lett.  A \vol 189 \yr 1994
 \pages  27-31 \endref

\ref \key Mat4  \bysame
       \paper MKdV Equation and  Submanifold Quantum Mechanics
      \jour  Sorushiron-Kenkyu \vol 94 \yr 1994 \pages
 A72-A75 \endref

\ref \key Mat5 \bysame
      \paper Anomaly on a Submanifold System:
       New Index Theorem related to
      a Submanifold System
      \jour  J.  Phys.  A.  \vol 28 \yr 1995 \pages
1399-1412 \endref

\ref \key Mat6 \bysame
       \paper The Physical Realization of the
 Jimbo-Miwa Theory of the
       Modified Korteweg-de Vries
  Equation on a Thin Elastic Rod: Fermionic Theory
        \jour Int.  J.  Mod.  Phys.  A \vol 10 \yr 1995
         \pages  3091-3107 \endref

\ref \key Mat7  \bysame
\jour Thesis of Tokyo Metroplitan University \yr 1996
 \paper On the Relation between Modified
  KdV Solitons and Dirac Fields on a
 Thin Elastic Rod \endref


\ref \key Mat8 \bysame
       \paper Constant Mean Curvature Surface and
Dirac Operator
       \jour J. Phys.  A \vol 30 \yr 1997
\pages  4019-4029 \endref


\ref \key Mat9 \bysame
      \paper Immersion Anomaly of Dirac Operator
on Surface in $\RR^3$
      \jour   Rev. Math. Phys. \yr 1999 \vol 2
\pages 171-186 \endref
\ref \key Mat10 \bysame
      \paper Dirac Operator of a Conformal
              Surface Immersed in $\RR^4$:
      Further Generalized Weierstrass
      Relation
      \jour  Rev. Math. Phys. \vol 12 \pages 431-444  \endref

\ref \key MT \by S.~Matsutani and H.~Tsuru
      \paper Physical relation between quantum mechanics and
       soliton on a thin elastic rod
       \jour Phys.  Rev.  A \vol 46 \yr 1992 \pages
          1144-1147 \endref


\ref \key P \by A. M. Polyakov \book Gauge Fields and Strings
\publ Harwood Academic Publ \publaddr London
 \yr 1987
\endref

\ref \key PP \by F.~Pedit and U.~Pinkall
\jour Doc. Math. J. DMV Extra Volume of ICM \vol II
\pages 389-400 \paper Quaternionic analysis on
Riemann surfaces and differential geometry  \yr 1998
\endref


\ref \key Tai1 \by  I.~A.~Taimanov
\paper Modified Novikov-Veselov equation and
differential geometry of surface \jour Trans. of
the Amer. Math. Soc.,
Ser.2 \vol 179 \yr 1997 \pages 133-151\endref


\ref \key Tai2 \bysame
\paper The Weierstrass representation of closed surfaces
in $\RR^3$ \jour Functional Anal. Appl. \yr 1998
      \vol 32\issue 4  \endref


\ref \key Tas \by  T.Tasaka
\book NijiKeishiki II (Quadratic Form II) \yr 1970
\publ Iwanami \publaddr Tokyo \lang japanese \endref

\ref \key Tr \by A. Trautman \jour  Acta Physica Polonica B
\paper Dirac operator on a hypersurface
\yr 1995 \pages 1283-1310
\endref

\endRefs
\bye